\documentclass[graybox]{svmult}

\usepackage{type1cm}        % activate if the above 3 fonts are
\usepackage{makeidx}         % allows index generation
\usepackage{graphicx}        % standard LaTeX graphics tool
\usepackage{multicol}        % used for the two-column index
\usepackage[bottom]{footmisc}% places footnotes at page bottom
\usepackage{amsmath}
\usepackage{hyperref}
\usepackage{cleveref}

\usepackage{tikz}
\usetikzlibrary{shapes,arrows,calendar,matrix,backgrounds,folding,calc,positioning,spy, colorbrewer}
\usetikzlibrary{positioning,shapes,arrows,calc, arrows.meta, patterns}

\usepackage{pgfplots}
\usepgfplotslibrary{colorbrewer}
\pgfplotsset{compat=1.15}

\pgfplotsset{
cycle list/Set1-5,
cycle multiindex* list={
mark list*\nextlist
Set1-5\nextlist
},
}

\usepackage{newtxtext}       % 
\usepackage[varvw]{newtxmath}       % selects Times Roman as basic font

\makeindex             % used for the subject index
\begin{document}

\makeatletter
\let\origsection\section
\renewcommand\section{\@ifstar{\starsection}{\nostarsection}}

\newcommand\nostarsection[1]
{\sectionprelude\origsection{#1}\sectionpostlude}

\newcommand\starsection[1]
{\sectionprelude\origsection*{#1}\sectionpostlude}

\let\origsubsection\subsection
\renewcommand\subsection{\@ifstar{\starsubsection}{\nostarsubsection}}

\newcommand\nostarsubsection[1]
{\sectionprelude\origsubsection{#1}\sectionpostlude}

\newcommand\starsubsection[1]
{\sectionprelude\origsubsection*{#1}\sectionpostlude}

\newcommand\sectionprelude{%
  \vspace{-2.25em}
}

\newcommand\sectionpostlude{%
  \vspace{-1.25em}
}
\makeatother

\title*{A Higher Order Unfitted Space-Time Finite Element Method for Coupled Surface-Bulk problems}
\titlerunning{A Higher Order Unfitted Space-Time FE Method for Coupled Surface-Bulk problems}
\author{Fabian Heimann\orcidID{0000-0002-8969-6504}}
\institute{Institute for Numerical and Applied Mathematics, University of Göttingen, \email{f.heimann@math.uni-goettingen.de}}
\maketitle

\abstract{We present a higher order space-time unfitted finite element method for convection-diffusion problems on coupled (surface and bulk) domains. In that way, we combine a method suggested by Heimann, Lehrenfeld, Preuß (SIAM J. Sci. Comput. 45(2), 2023, B139 - B165) for the bulk case with a method suggested by Sass, Reusken (Comput. Math. Appl. 146(15), 2023, 253-270) for the surface case. The geometry is allowed to change with time, and the higher order discrete approximation of this geometry is ensured by a time-dependent isoparametric mapping. The space-time discretisation approach allows for straightforward handling of arbitrary high orders. In that way, we also generalise results of Hansbo, Larson, Zahedi (Comput. Methods Appl. Mech. Engrg. 307, 2016, 96-116) to higher orders. The convergence of the proposed higher order discretisations is confirmed numerically.}
\section{Introduction}
Finite Element methods are well-established tools for the simulation of physical phenomena on complex geometries. In this article, we focus on coupled surface-bulk problems on moving domains and develop a novel finite element method to solve an example equation from this category, the convection-diffusion equation. More generally, coupled dimensional convection-diffusion-type systems have been used to predict the dynamics of biological cells, and more generally applied in multi-phase flows involving surfactants. These applications motivate to develop methods which are robust under significant changes of potentially complex geometries. In the last decades, several such methods were developed under the name of \textit{Unfitted Finite Element} methods.\cite{burman15} The fundamental idea of these methods is the use of a fixed background mesh, which delivers basis functions for the discrete method. Challenges arise in relation to numerical integration and stabilisation. Two subcategories of the Unfitted FE methods are surface and bulk methods, where in the latter case the spatial dimensions of domain of interest and mesh align, whilst in the surface case the background mesh is of one dimension higher. The specific method we put forward in this article can be regarded as the combination of the method \cite{HLP2022} developed specifically for a convection-diffusion in the bulk and the method \cite{sassreusken23} modelling a corresponding surface equation.

One particular challenge of Unfitted Finite Element methods in that of numerical integration of higher order accuracy on implicitly defined domains. To be specific, we employ the levelset framework and assume that the implict geometry is given in terms of a levelset function $\phi$, $\Omega = \{\phi <0 \}$. To perform numerical integration of higher order accuracy on such domains, we use the technique of an isoparametric mapping $\Theta_h$. It involves---in the case of a simplicial spatial mesh---a linear reference geometry, on which numerical integration can be solved by tessellation, which is then mapped by $\Theta_h$ in order to arrive at a high order accurate discrete geometry.\footnote{In the case of quadrilateral meshes, a slightly modified version of the construction applies, c.f. \cite{HL_ENUMATH_2019} for more details. Throughout this text, we consider simplicial meshes.}

We aim at solving the full time-dependent convection-diffusion equation, and opt for a space-time method in order to solve for a time dependent solution. The space-time finite elements of consideration will have the structure $T \times I_n$, where $T \in \mathcal{T}_h$ would denote a usual spatially Unfitted Finite Element, and $I_n = [t_{n-1}, t_n]$ the time interval, so that $\Delta t = t_n - t_{n-1}$. These methods have the benefit that they yield a uniform formulation as it regards to the discretisation order in time, so that our approach will naturally leed to a class of methods of arbitrary high order in space and time.

Apart from the references given so far, we would like to mention the literature \cite{HLZ16, Z18}, where similar methods have been developed involving a slightly different stabilisation method and a smaller set of investigated discrete orders. Also, in \cite{massing18} a coupled surface-bulk problem for the stationary case has been investigated.

The remainder of this paper is structured as follows: In the following short \Cref{sec:model_prob}, we introduce the model problem of consideration. Afterwards, in \Cref{sec:discr_method}, we introduce the discrete method to solve this problem. Then, in \Cref{sec:num_exp}, we give numerical results for example problems. Finally, we give a short summary and an outlook in \cref{sec:summ_outlook}.
\section{Model Problem} \label{sec:model_prob}
\begin{figure}[tb] \begin{center}
    \begin{tikzpicture}[rotate=-90]
  \draw (0,0) rectangle (4,5);
  \draw (2,2) circle (1cm);
  \node at (2,2) {$\Omega_2(t)$};
  \node at (2,4) {$\Omega_1(t)$};
  \draw[Set1-B, ->] (3,2) to[bend right=30] (3.3,2.3);
  \node[Set1-B] at (3.3,1.9) {$w$};
  \node at (1.5,0.9) {$\Gamma$};
  \node at (0.3,2.5) {$\tilde \Omega = \Omega_1(t) \cup \Omega_2(t)$};
 \end{tikzpicture}
 \end{center} \vspace{-0.4cm}
 \caption{Illustration of the decomposition of the fixed domain $\tilde \Omega$ into $\Omega_1(t)$ and $\Omega_2(t)$. For the surface $\Gamma(t) = \partial \Omega_2(t)$.}
 \label{fig:dom_illu}
 \vspace{-0.5cm}
\end{figure}
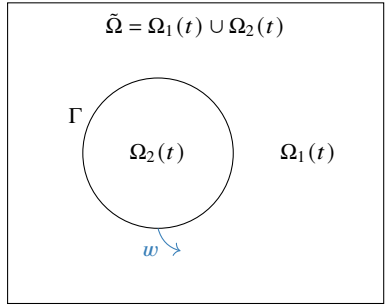

To introduce the model problem, we start with setting up some geometry conventions, following our initial motivational examples of a biological cell or a bubble in a multi-phase flow. (C.f. \Cref{fig:dom_illu} below for a sketch.) Let the shape of this interior object be denoted by a smooth $\Omega_2(t) \subseteq \mathbb{R}^d$, for some time $t$ of interest in an appropriate time interval, $t \in [0,T]$. The convection-diffusion surface problem should operate on the boundary of this domain, motivating the following definition for the space-time boundary:
\begin{equation}
 \Gamma = \bigcup_{t\in [0,T]} \partial \Omega_2(t) \times \{t\}, \quad \Gamma(t) = \partial \Omega_2(t).
\end{equation}
Apart from the surface convection-diffusion problem, we also want to include a convection-diffusion problem in the surrounding bulk domain in the model. Hence, we assume that there is an also smooth domain of interest $\Omega_1(t) \subseteq \mathbb{R}^d$, given such that $\Omega_1(t) \cup \Omega_2(t) = \tilde \Omega$ for some fixed connected polygonal background domain $\tilde \Omega$. We also write $\Omega(t) = \Omega_1(t)$. This motivates the following notion of a space-time domain.
\begin{equation}
 Q = \bigcup_{t\in [0,T]} \Omega(t) \times \{t\}.
\end{equation}
The concentration of the relevant species on the surface now shall be denoted by a function $u_S=u_S(x,t)$, whereas $u_B = u_B(x,t)$ will be used for the bulk concentration function. In the bulk, $\Delta$ and $\nabla$ are assumed to be the usual cartesian differential operators inherited from $\mathbb{R}^d$. For the surface, let $\vec{n}$ denote the unique outward pointing normal to $\Omega_2(t)$ on $\Gamma(t)$. Then, we introduce a projection operator $\vec{P}_\Gamma$ and the corresponding surface differential operators as follows:
\begin{align}
 \vec{P}_\Gamma &= \vec{I} - \vec{n} \otimes \vec{n}, \quad && \nabla_\Gamma = \vec{P}_\Gamma \cdot \nabla\\
 \operatorname{div}_\Gamma \vec{v}& = \mathrm{tr}(\vec{v} \otimes \nabla_\Gamma), \quad &&\Delta_\Gamma = \operatorname{div}_\Gamma(\nabla_\Gamma \dots)
\end{align}
In terms of these definitions, the coupled surface-bulk convection diffusion problem in strong form reads as follows: Find $u_B\colon Q \to \mathbb{R}$ and $u_S \colon \Gamma \to \mathbb{R}$ s.t.
  \begin{align}
 \partial_t u_B + \vec{w} \cdot \nabla u_B - k_B \Delta u_B &= 0 && \textnormal{in } \Omega(t) ~ \forall t \in (0,T], \nonumber\\
 - k_B \nabla u_B \cdot \vec{n}_{\partial \tilde \Omega} &= 0 && \textnormal{on } \partial \tilde \Omega ~ \forall t \in (0,T], \nonumber\\
 - k_B \nabla u_B \cdot \vec{n} &= f_{\text{coupl}} && \textnormal{on } \Gamma(t) ~ \forall t \in (0,T], \label{strongformproblem}\\
 \partial_t u_S + \vec{w} \cdot \nabla u_S + (\operatorname{div}_\Gamma \vec{w}) u_S - k_S \Delta_\Gamma u_S &= f_{\text{coupl}} && \textnormal{on } \Gamma(t) ~ \forall t \in (0,T], \nonumber \\
 u_B(\cdot, 0) = u_{B,0} \textnormal{ in } \Omega(0), \quad u_S(\cdot, 0) &= u_{S,0} &&\textnormal{ on } \Gamma(0).\nonumber
\end{align}
for suitable initial data $u_{B,0}$, $u_{S,0}$, coupling force right-hand side $f_{\text{coupl}}$, diffusion constants $k_B$, $k_S$, and a convection velocity field $\vec{w}$ with $\mathrm{div}(\vec{w}) = 0$. Also, we assume that the movement of $\Gamma(t)$ matches with $\vec{w}$.

For the coupling, we distinguish the general non-linear ansatz $f_{\text{coupl}} = b_B u_B - b_S u_S - b_{BS} u_B u_S$ with appropriate constants $b_B, b_s$, and $b_{BS}$ from the linear ansatz, where only the constants $b_S$ and $b_B$ are to be chosen as above, and $b_{BS} = 0$. The linear case is---up to  naming of constants---called the Henry model, whereas the non-linear ansatz represents the Langmuir model. For more details on the motivation and the scope of application of this framework, we refer to the literature. (See e.g. \cite{HLZ16} and the references therein, \cite{grossetal15} which also contains a mathematical well-posedness analysis of a related coupled problem.)
\section{Discrete Method} \label{sec:discr_method}
We now derive a discrete method for this strong form problem. In line with the unfitted paradigm, we assume that on our time-independent background domain $\tilde \Omega$ a family of simplicial meshes $\mathcal{T}_h$ of decreasing mesh size $h$ is given. The time interval of interest, $[0,T]$, is assumed to be subdivided into time intervals or steps $I_n = [t_{n-1},t_n], n =1,\dots,N, t_0 = 0, t_N = T$ of equal length $\Delta t = t_n - t_{n-1}$. With the space-time method, we intend to solve problems on each of these time steps $I_n$ in order to combine the flexibility of chosing a discretisation order with only moderately sized discrete linear algebra systems.
\subsection{Geometry handling}
In order to obtain a higher order accurate and computationally feasible discrete geometry, we employ the time-dependent variant of an isoparametric mapping. \cite{HLP2022, L_CMAME_2016, heimann2023geometrically} Fundamentally, it is assumed that $\Omega_2$ is given in accordance with a levelset function $\phi(x,t)$, i.e.
\begin{equation}
 Q = \{ (x,t) \in \tilde \Omega \times [0,T] \, | \, \phi(x,t) <0 \},
\end{equation}
and we are given a $\mathcal{T}_h$-piecewise linear-in-space and temporally higher order interpolation thereof, $\phi^{\text{lin}}$. The induced linear reference configuration, which we denote as $Q^{\text{lin}}$ / $\Omega^{\text{lin}}(t)$ is instrumental for the numerical integration procedure as for fixed times of interest, $\Omega^{\text{lin}}(t)$ is a polygonal. This enables us to construct an algorithm of numerical integration on $Q^{\text{lin}}$, which is described in detail in \cite{HLP2022}. Despite this computational advantage, $Q^{\text{lin}}$ is only a spatially second order accurate approximation of $Q$. To improve upon this, the domain is mapped by a function $\Theta_h^{\text{st}}\colon \tilde \Omega \times [0,T] \to \tilde \Omega \times [0,T]$, which preserves the time coordinate. The relevant discrete domains are defined as images under this mapping:
\begin{align}
  \Omega^h(t) &= \Theta_h^{\text{st}} (\Omega^{\text{lin}}{t}, t),
  &Q^{\text{lin}, n} &= {\bigcup}_{t \in I_n} \Omega^{\text{lin}}(t) \times \{ t \},
%  Q^{\text{lin}} &= {\bigcup}_{n = 1}^{N} \Qlinn,
  &Q^{h,n} &= \Theta_h^{\text{st}}(Q^{\text{lin},n}) \label{eq:introQhn},\\
 \Gamma^h(t) &= \Theta_h^{\text{st}} (\partial \Omega^{\text{lin}}(t), t), \hspace{-0.1cm}
  &\Gamma^{\text{lin}, n} &= {\bigcup}_{t \in I_n} \partial \Omega^{\text{lin}}(t) \times \{ t \}, \hspace{-0.1cm}
%  Q^{\text{lin}} &= {\bigcup}_{n = 1}^{N} \Qlinn,
 &\Gamma^{h,n} &= \Theta_h^{\text{st}}(\Gamma^{\text{lin}, n}).
\end{align}
We refer the reader to \cite{HLP2022} for more details on the computational aspects, as well as to \cite{heimann2023geometrically} for a mathematical analysis of the higher order accuracy of the so-defined computational domains.
\subsection{Discrete spaces}
For the discrete spaces, we decide to introduce a notion of active elements within each time step $I_n$. In this way, only the respective active elements will induce degrees of freedom in the linear system. We define for bulk and surface respectively
 \begin{align*}
 \mathcal{E}(Q^{\text{lin}, n}) :&= \bigcup \{ T \times I_n \, | \, (T \times I_n) \cap Q^{\text{lin}, n} \neq \varnothing, T \in \mathcal{T}_h \}, \\
  \mathcal{E}(\Gamma^{\text{lin}, n}) :&= \bigcup \{ T \times I_n \, | \, (T \times I_n) \cap \Gamma^{\text{lin},n} \neq \varnothing, T \in \mathcal{T}_h \}.
 \end{align*}
 On these regions, which are illustraed in \Cref{fig:discr_reg}, we introduce our space-time finite element spaces of order $k_s$ in space and $k_t$ in time:
\begin{align*}
  W_h^{n,b} &= \{ v \in V_h^{k_s} \otimes \mathcal{P}^{k_t}(I_n) \, | \, v \textnormal{ vanishes outs. of } \mathcal{E}(Q^{\text{lin}, n}) \textnormal { as discrete fn.} \} \circ (\Theta_h^{\text{st}})^{-1}, \\
  W_h^{n,s} &= \{ v \in V_h^{k_s} \otimes \mathcal{P}^{k_t}(I_n) \, | \, v \textnormal{ vanishes outs. of } \mathcal{E}(\Gamma^{\text{lin}, n}) \textnormal { as discrete fn.} \} \circ (\Theta_h^{\text{st}})^{-1}.
  \end{align*}
  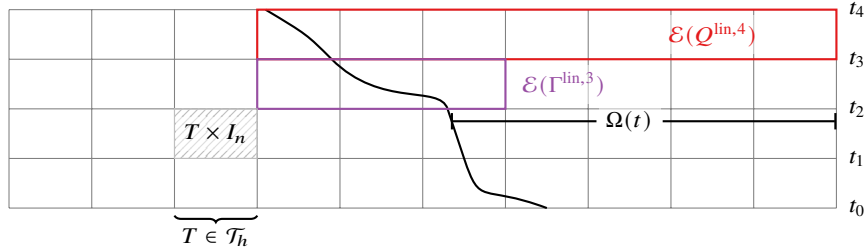
\begin{figure}  \vspace{-0.5cm} \centering
   \begin{tikzpicture}[xscale=1.1, yscale=0.66]
  \draw[gray] (0,0) grid[xstep=1, ystep=1] (10,4);
  \begin{scope}[xscale=-1, xshift=-10cm]
  \draw[thick] (3.5,0) .. controls (4.3, 0.5) and (4.3,0) .. (4.5,1)
                .. controls (4.6,1.5) .. (4.7,2)
                .. controls (4.8,2.5) and (5.5,2) .. (6.1,3)
                .. controls (6.4,3.5) .. (6.9,4);
   \end{scope}
   \draw [
    thick,
    decoration={
        brace,
        mirror,
        raise=0.1cm
    },
    decorate
] (2,0) -- (3,0)
node [pos=0.5,anchor=north,yshift=-0.15cm] {$T \in \mathcal{T}_h$};
%    \draw[gray!50!white, pattern = north east lines, pattern color = gray!50!white] (5,0) rectangle (7,1) -- (5,1) rectangle (6,2) -- (3,2) rectangle (6,3) -- (3,3) rectangle (4,4);
%    \draw[dash dot, Set1-B, thick] (7,0) -- (7,1) -- (6,1) -- (6,3) -- (4,3) -- (4,4) -- (10,4) -- (10,0) -- cycle;
%    \only<2>{\draw[very thick, dashed, Set1-C] (7,0) -- (7,1) (6,0) -- (6,3) (5,2) -- (5,3) (4,2) -- (4,4);
%    \draw[very thick, dashed, Set1-C!50!black] (8,0) -- (8,1) (7,2) -- (7,3) (8,2) -- (8,3);}

%    \draw[{Bar[]}-, thick] (5.0,1.15) --(6.8,1.15);
%    \draw[-{Bar[]}, thick] (8.1,1.15) -- (10, 1.15);
%    \node[fill=white]  at (7.5,1.15) {$\mathcal{E}(\Omega^{\text{lin},2})$};

%   \fill[blue!50!white, fill opacity=0.1] (1,4) rectangle (2,5);
   \node[Set1-A,fill=white, fill opacity=0.5, text opacity=1] at (8.5,3.5) {$\mathcal{E}(Q^{\text{lin}, 4})$};
%    \node[Set1-B,fill=white] at (8.5,3.6) {$\bigcup_{n=1}^N \mathcal{I}(Q^{\text{lin}, n})$};
%    \node[gray,fill=white] at (3,1.5) {$\bigcup_{n=1}^N \mathcal{E}(Q^{\text{lin}, n}) \backslash \mathcal{I}(Q^{\text{lin},n})$};
%    \node[Set1-C,fill=white] at (6.35,2.5) {$\mathcal{F}^{n}_R$};
%    \node[Set1-C!50!black] at (8.85,2.51) {$\mathcal{F}^{n, \text{ext}}_R \backslash \mathcal{F}^{n}_R$};

   \draw[Set1-A, thick] (3,3) -- (3,4) -- (10,4) -- (10,3) -- cycle;

   \node[Set1-D,fill=white, fill opacity=0.5, text opacity=1] at (6.7,2.5) {$\mathcal{E}(\Gamma^{\text{lin}, 3})$};
   \draw[Set1-D, thick] (3,2) -- (3,3) -- (6,3) -- (6,2) -- cycle;

   \draw[gray!50!white, pattern = north east lines, pattern color = gray!50!white] (2,1) rectangle (3,2);
   \node[fill=white, fill opacity=0.5, text opacity=1] at (2.5,1.5) {$T \times I_n$};

   \node at (10.25, 0) {$t_0$};
   \node at (10.25, 1) {$t_1$};
   \node at (10.25, 2) {$t_2$};
   \node at (10.25, 3) {$t_3$};
   \node at (10.25, 4) {$t_4$};

%    \node at (10.1, 1.75) {$t$};
   \draw[{Bar[]}-, thick] (5.34,1.75) --(7.1,1.75);
   \draw[-{Bar[]}, thick] (7.9,1.75) -- (10, 1.75);
   \node[fill=white, fill opacity=0.5, text opacity=1] at (7.5,1.75) {$\Omega(t)$};

 \end{tikzpicture}\vspace{-0.2cm}
 \caption{Illustration of the discrete regions $\mathcal{E}(Q^{\text{lin}, n})$, $\mathcal{E}(\Gamma^{\text{lin}, n}) $ in a one-dimensional example.}
 \label{fig:discr_reg}
  \vspace{-0.4cm}
  \end{figure}
  Here, $V_h^{k_s}$ denotes the usual continuous finite element space of order $k_s$ on the mesh $\mathcal{T}_h$ and $\mathcal{P}^{k_t}(I_n)$ the space of scalar one-dimensional polynomials on $I_n$ up to order $k_t$. Finally, these spaces are taken together for the overall discrete space $W_h^n := W_h^{n,b} \otimes W_h^{n,s}$.
\subsection{Bilinear- and linear form for the linear case}
 A weak form of the problem \cref{strongformproblem} is obtained by multiplying with a test function scaled with the material constant ($b_B v_B$ and $b_S v_S$ for bulk and surface respectively) and integrating over the space-time domain of interest of the local time slice, followed by an application of partial integration. Moreover, we are interested in the discrete variational form, where the function spaces and integration domains are replaced with appropriate discrete counterparts. This yields:
 Find $u=(u_B, u_S) \in W_h^n$ such that $\forall v=(v_S, v_B) \in W_h^n$
\begin{equation}
  b_B B^n_b(u_B,v_B) + b_S B^n_s(u_S,v_S) + B_{\text{upw}}^n(u,v) + J^n(u,v) + B_{\text{coup}}(u,v) = f^n_{\text{upw}}(v), %\label{discreteproblem}
\end{equation}
where
 \begin{align*}
   B^n_b(u_B,v_B) \! & :=\! (\partial_t u_B \!+\! \vec{w} \cdot \nabla u_B, v_B)_{Q^{h,n}}
   \!+ k_B (\nabla u_B, \nabla v_B)_{Q^{h,n}},
   \end{align*}
   \begin{align*}
   B^n_s(u_S,v_S) \! & :=\! (\partial_t u_S \!+\! \vec{w} \cdot \nabla u_S + u_S \operatorname{div}_\Gamma \vec{w}, v_S)_{\Gamma^{h,n}}
   \!+ k_S (\nabla_\Gamma u_S, \nabla_\Gamma v_S)_{\Gamma^{h,n}}, \\
   \hspace*{-0.25cm} B_{\text{upw}}^{n}(u_B,v_B)\! &:= \! b_B ({u}^{n-1}_{B,+}\!, v^{n-1}_{B, +})_{\Omega^{h}\!(t_{n-1}\!)} + b_S ({u}^{n-1}_{S,+}\!, v^{n-1}_{S, +})_{\Gamma^{h}\!(t_{n-1}\!)},\\
   f_{\text{upw}}^{n}(v)\! &:= \! b_B ( \Pi_u {u}^{n-1}_{B,-}\!, v^{n-1}_{B,+})_{\Omega^{h}\!(t_{n-1}\!)} + b_S ( \Pi_u {u}^{n-1}_{S,-}\!, v^{n-1}_{S,+})_{\Gamma^{h}\!(t_{n-1}\!)} \\
   B_{\text{coup}}(u,v) &= (b_B u_B - b_S u_S, b_B v_B - b_S v_S)_{\Gamma^{h,n}}
 \end{align*}
%\end{subequations}
with $u^n_{+} := \lim_{s \searrow t_n} u(\cdot, s)$, $u^n_{-} := \lim_{s \nearrow t_n} u(\cdot, s)$, and $(\cdot, \cdot)_S$ denoting the usual $L^2$ inner product on the geometric domain $S$. Note that we used an upwind (in time) stabilisation in order to enforce continuity along time slice boundaries in a weak sense. The (coupling) bilinear form is a bilinear form as $b_{BS} = 0$ is to be considered in the linear case. With $\Pi_u$ we denote a shifted evaluation discrete mapping, which is necessary for potentially discontinuous discrete regions. It effectively projects the old discrete function $u^n_-$ from some space $V_h^{k_s} \circ (\Theta_h^{\text{st}}(\cdot,t_n)^-)$ into a potentially slightly differing $V_h^{k_s} \circ (\Theta_h^{\text{st}}(\cdot,t_n)^+)$. (C.f. \cite{HLP2022, LL_ARXIV_2021, heimann2023geometrically} for details.) With $J^n$, we denote the direct Ghost penalty stabilisation introduced in \cite{HLP2022} plus the normal gradient stabilisation known from \cite{sassreusken23}.
\begin{align}
 J^n(u, v) = &\sum_{F \in \mathcal{F}_n} \int_{I_n} \int_{\omega_F} \frac{\gamma_B}{h^2} \left( 1 + \frac{\Delta t}{h} \right) [u_B]_{\omega_F} [v_B]_{\omega_F} \nonumber \\
 & + \int_{\Theta_h^{\text{st}}(\mathcal{E}(\Gamma^{\text{lin},n}))} \gamma_S (\vec{n_h} \cdot \nabla u_S) (\vec{n_h} \cdot \nabla v_S),
\end{align}
where $\mathcal{F}_n$ is an appropriate facet set containing at least the inner spatial facets of $\mathcal{E}(\Gamma^{\text{lin},n})$, $\omega_F$ denotes the union of the elements which contain the complete facet $F$ (mapped with the isoparametric mapping), and $[ \dots ]_{\omega_F}$ the jump between a function and the polynomial extension of its counterpart from the other element in $\omega_F$. (C.f. \cite{HLP2022} for more details.) For the normal gradient stabilisation, $\vec{n_h}$ denotes the discrete approximate normal vector, i.e. the normal of $\Gamma^{\text{lin},n}(t)$ transformed with $\Theta_h^{\text{st}}(\cdot,t)$. (C.f. \cite{sassreusken23, grande2018analysis})
Overall, the stabilisation ensures well-posedness of the linear system in the presence of bad cut configurations, which are configurations where the implicitly constructed elements (cut parts of shape-regular $T \in \mathcal{T}_h$) do not feature the inverse inequalities needed for the numerical analysis. Then, the discrete control can be obtained from interior elements by walking along facets from $\mathcal{F}_n$, c.f. \cite{burman2010ghost, burman15, preuss18, heimann20}. The normal gradient stabilisation is known from respective surface problems. \cite{grande2018analysis, burman18cutfemcodim, sassreusken23, reusken2024analysis}
\subsection{Bilinear- and linear form for the non-linear case}
For the nonlinear case, we linearise the problem in order to obtain the following formulation. We define a residual semi-linear form $F(u,v)$ for $u,v \in W_h^n$ as follows, where the argument $v$ shall be understood as the linear form argument and $u$ is a current fixed iteration candidate function:
\begin{align}
 F(u,v)= &b_B B_b^n(u_B,v_B) + b_S B^n_s(u_S, v_S) + B^n_{\text{upw}}(u,v)-f^n_{\text{upw}}(v) +  J^n(u,v) \nonumber \\
 &+B_{\text{coup}}(u,v) + b_{BS} (u_B u_S, b_S v_S - b_B v_B)_{\Gamma^{h,n}}.
\end{align}
The reader easily confirms that $F(u,v) = 0$ if $u$ was the solution to \cref{strongformproblem}, involving arbitrary $b_{BS}$. Now, in order to obtain this root, we define the bilinear form $DF(w,v)$ as follows:
\begin{align}
 DF(w,v) = &b_B B_b^n(w_B,v_B) + b_B B^n_s(w_S, v_S) + B^n_{\text{upw}}(w,v) + J^n(w,v) +B_{\text{coup}}(u,v) \nonumber \\
 &+ b_{BS} (w_B u_S + u_B w_S, b_S v_S - b_B v_B)_{\Gamma^{h,n}}.
\end{align}
The numerical method for determining $u$ then proceeds as follows: As the initial guess, we assume the generally time-dependent function $u$ to be the fixed  initial data $\Pi_u u = (\Pi_u u_S, \Pi_u u_B)$ (either from initial data or previous timestep). We assemble $F$ and $DF$ according to the above formulas, and calculate $w = (\mathcal{D} F)^{-1} \mathcal{F}$, where by $\mathcal{D} F$ we denote the finite element matrix associated to $DF(\cdot, \cdot)$ by the Galerkin isomorphism, and by $\mathcal{F}$ the associated vector to $F(u,\cdot)$. $w$ is subtracted from the discrete vector of $u$, and the procedure is repeated until the $L^2$ norm of $w$ (as a discrete vector) is below $10^{-9}$, which indicates convergence. This Newton method is not guaranteed to converge and we will see in numerical experiments that specifically in the context of coarse curved meshes and large time steps, challenges arise.
\section{Numerical Examples} \label{sec:num_exp}
\begin{figure}
 \begin{center}
  \includegraphics[width=0.49\textwidth]{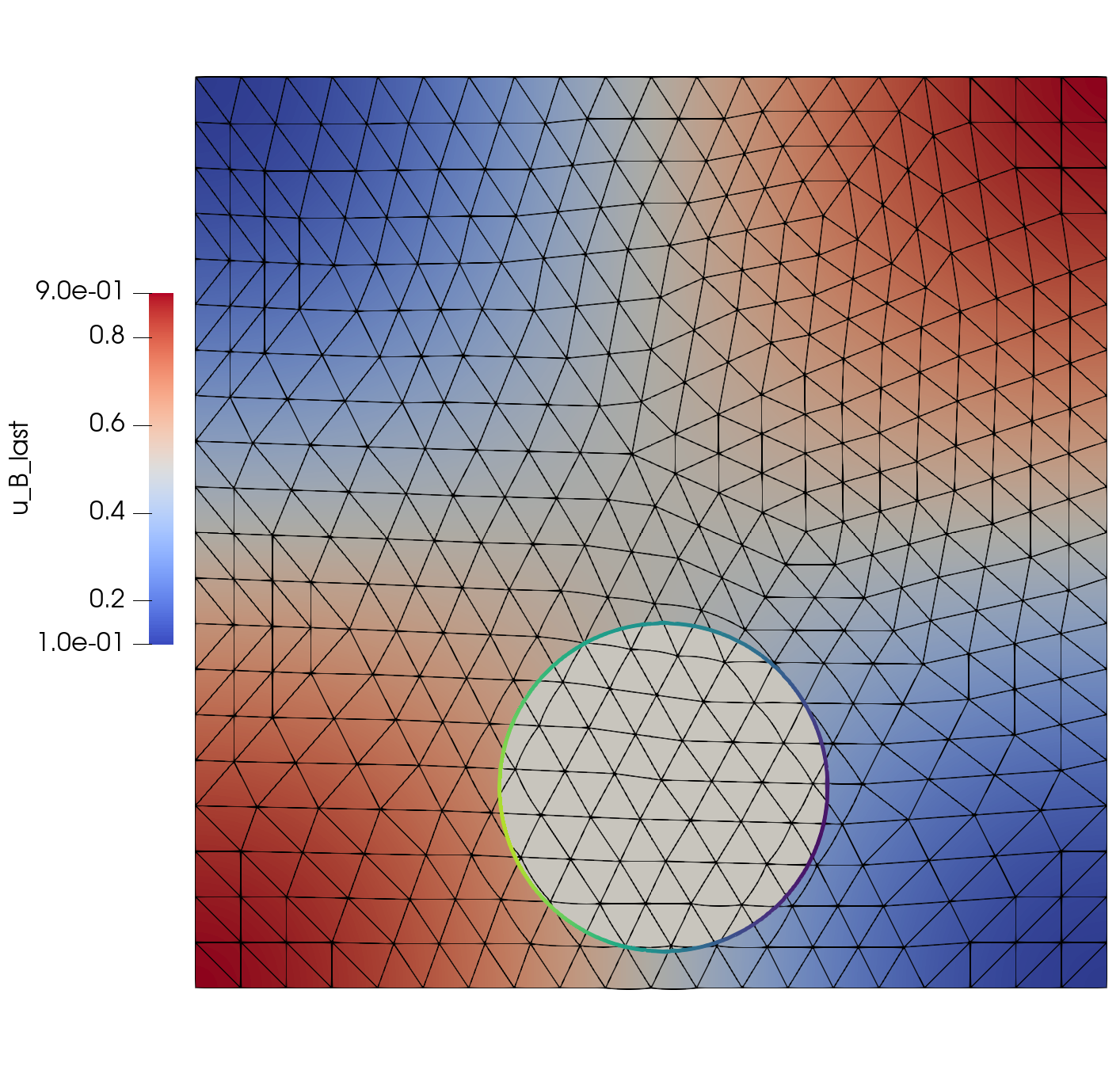}
  \includegraphics[width=0.49\textwidth]{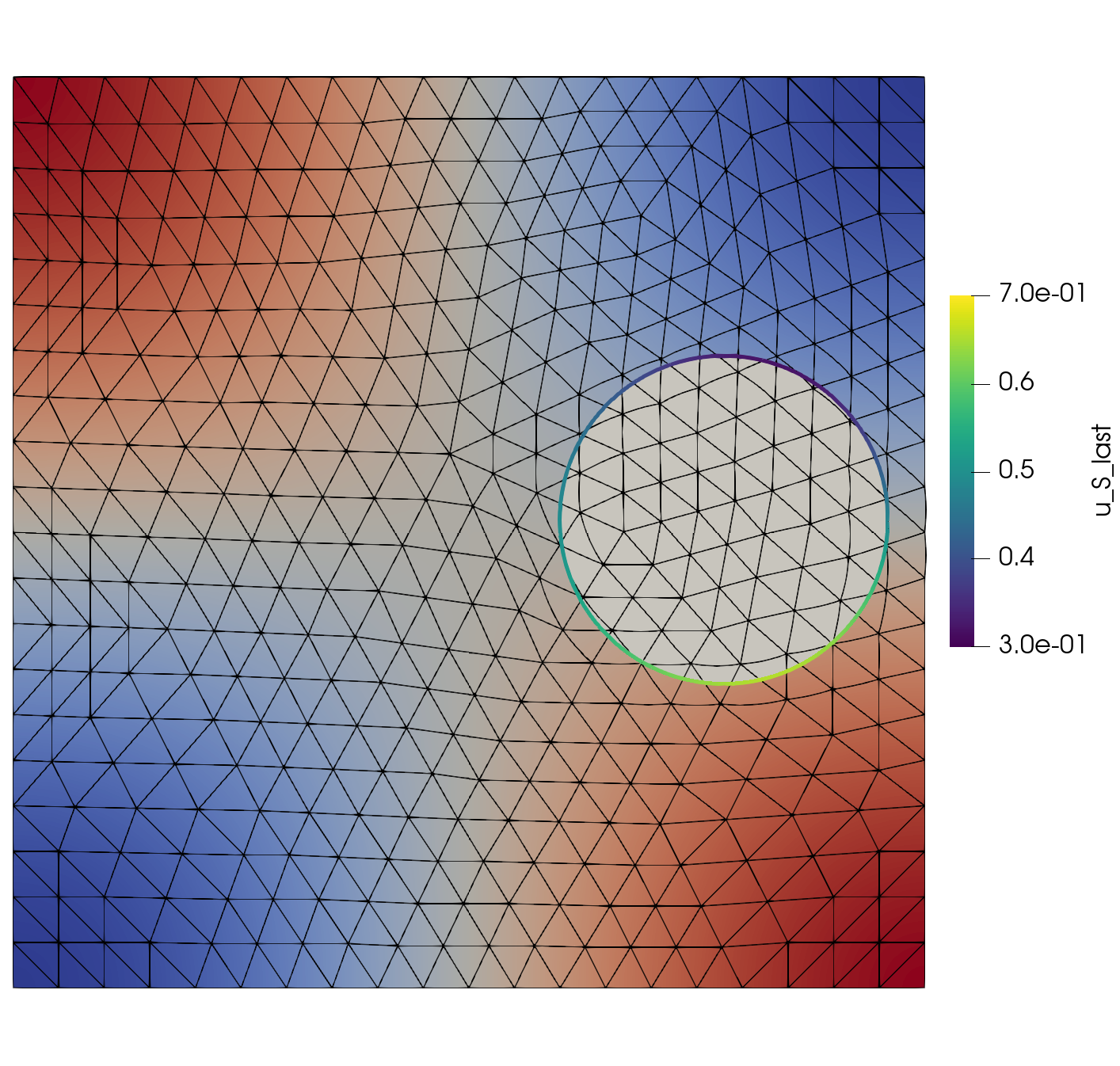}
 \end{center}\vspace{-0.5cm}
 \caption{Plots of the manufactured solution $u_B$ on the bulk, and $u_S$ on the surface. They respectively have their colour bar, indicating some coupling between higher/ lower concentration values. \textit{Left}: $u_B(x,0)$, $u_S(x,0)$. \textit{Right}: $u_B(x,T=0.5)$, $u_S(x,T=0.5)$.}
 \label{snapshots_fun}
\end{figure}

We investigate numerically the high-order convergence behaviour of the suggested method. To this end, we fix a two-dimensional test geometry of a square, $\tilde \Omega = [0,1]^2$, with a circular hole moving around with a convection field $\vec{w}$ in the time interval $[0, T= 0.5]$. In particular
 $\phi(x,t) = 0.18 - \sqrt{(x-x_c(x,t))^2 + (y-y_c(x,t))^2}, \vec{w} = ( \pi (0.5 - y), \pi (x-0.5))^T$,
 $x_c (x,t) = 0.5 + 0.28 \sin (\pi t), ~~~ y_c (x,t) = 0.5 + 0.28 \cos (\pi t)$.
For the diffusion and coupling constants, we choose $k_B = 0.01$, $k_S = 1$, $b_B = b_S = 1$, ($b_{BS} = 1$ in the non-linear case), $\gamma=\gamma_S=\gamma_B=0.05$ and the right hand-sides are chosen in accordance with the manufactured solution $u_B = 0.5 + 0.4 \cos (\pi x) \sin (\pi y) \cdot \cos(2 \pi t)$. Following the manufactured solution paradigm, we calculate from this $f_{\text{coupl}} = -k_B \nabla u_B \cdot \mathbf{n}$, $u_S = (b_B u_B - f_{\text{coupl}})/(b_S + b_{BS} u_B)$. $u_B$ and $u_S$ are displayed in \Cref{snapshots_fun} for illustration.

We implement the discrete method in our software \texttt{ngsxfem} \cite{xfem_joss}, an Add-On to NGSolve for Unfitted Finite Elements. In particular, it contains routines for numerical integration and the higher order space-time isoparametric mapping. Reproduction data is available at \url{https://gitlab.gwdg.de/fabian.heimann/repro-ho-unf-space-time-coupled}. We investigate the numerical error on a series of simultaneous refinements in space and time, i.e. $h = 0.2 \cdot 0.5^i, \Delta t = 0.5^{i+2}$. Moreover, we fix an order parameter $k = k_s = k_t$ to match the most straightforward choice of order combinations already known from the bulk problem. $k$ is also the parameter of geometrical accuracy in space and time for the isoparametric mapping ($q_s= q_t = k$ following the notation of \cite{HLP2022}). We measure the error respectively in an $L^2$ norm at the final time on bulk and surface: $\| u_B - u_{B,h} \|_{L^2(\Omega(T))}^2 = (u_B - u_{B,h}, u_B - u_{B,h})_{\Omega^h(T)}$, $\| u_S - u_{S,h} \|_{L^2(\Gamma(T))}^2 = (u_S - u_{S,h}, u_S - u_{S,h})_{\Gamma^h(T)}$.

\begin{figure}
 \begin{tikzpicture}[scale=0.7]
    \begin{semilogyaxis}[ xlabel=$i$, ylabel={numerical error $L^2(\Omega(T))/ L^2(\Gamma(T))$}, legend entries ={ $k=1$ $\Omega$, $k=2$ $\Omega$, $k=3$ $\Omega$, $k=4$ $\Omega$, $k=1$ $\Gamma$, $k=2$ $\Gamma$, $k=3$ $\Gamma$, $k=4$ $\Gamma$
%      {[text width=2.5cm, text depth=] $O(h^{k+1})= O(\Delta t^{k+1})$}
    }, legend style={at={(0.5,-0.1)},anchor=north,legend columns=4, draw=none, legend cell align=left, fill=none},
    x label style={at={(axis description cs:0.5,0.05)},anchor=center},
%     ymax=3, ymin=3e-11
     ]
     \addplot table {out/conv_lin_ks1_kt1_ktls1_from0_to7.dat};
     \addplot table {out/conv_lin_ks2_kt2_ktls2_from0_to6.dat};
     \addplot table {out/conv_lin_ks3_kt3_ktls3_from0_to5.dat};
     \addplot table {out/conv_lin_ks4_kt4_ktls4_from0_to4.dat};
     \pgfplotsset{cycle list shift=1}
     \addplot table[x index = 0, y index =2] {out/conv_lin_ks1_kt1_ktls1_from0_to7.dat};
     \addplot table[x index = 0, y index =2] {out/conv_lin_ks2_kt2_ktls2_from0_to6.dat};
     \addplot table[x index = 0, y index =2] {out/conv_lin_ks3_kt3_ktls3_from0_to5.dat};
     \addplot table[x index = 0, y index =2] {out/conv_lin_ks4_kt4_ktls4_from0_to4.dat};
     \addplot[gray, dashed, domain=0:7] {(1/2^(x+2.5)))^2};
     \addplot[gray, dashed, domain=0:6] {(1/2^(x+3)))^3};
     \addplot[gray, dashed, domain=0:5] {(1/2^(x+2.5)))^4};
     \addplot[gray, dashed, domain=0:4] {(1/2^(x+2.5)))^5};
    \end{semilogyaxis}
     \node[scale=0.75] at (4.5,5.25) {$O(h^{k+1})= O(\Delta t^{k+1})$};
     \draw[scale=0.75, gray, dash=on 2.25pt off 2.25pt phase 0pt, line width=0.4*0.75pt] (6.25/0.75,5.25/0.75) -- (6.8/0.75,5.25/0.75);
   \end{tikzpicture} %\vspace*{-1cm}
        \begin{tikzpicture}[scale=0.7]
    \begin{loglogaxis}[ xlabel={runtime [s]}, ylabel=numerical error $L^2(\Omega(T))$, legend entries ={ $k=1$, $k=2$, $k=3$, $k=4$
%      {[text width=2.5cm, text depth=] $O(h^{k+1})= O(\Delta t^{k+1})$}
    }, legend style={anchor=north,legend columns=1, draw=none, legend cell align=left, fill=none}, legend pos =south west,
%     x label style={at={(axis description cs:0.5,0.05)},anchor=center},
%     ymax=3, ymin=3e-11
     ]
     \addplot table[x index = 5, y index =1] {out/conv_lin_ks1_kt1_ktls1_from0_to7.dat};
     \addplot table[x index = 5, y index =1] {out/conv_lin_ks2_kt2_ktls2_from0_to6.dat};
     \addplot table[x index = 5, y index =1] {out/conv_lin_ks3_kt3_ktls3_from0_to5.dat};
     \addplot table[x index = 5, y index =1] {out/conv_lin_ks4_kt4_ktls4_from0_to4.dat};
    \end{loglogaxis}
    \draw[opacity=0] (0,0) -- (0,-1.5);
   \end{tikzpicture}
   \caption{Numerical convergence results for the linear model system ($b_{BS} = 0$) and different polynomial orders $k=k_s=k_t=q_s=q_t$. Refinement by $i$ refines both space and time simultaneously. \textit{Left}: The discrete error is measured in bulk and surface norms vs. refinements $i$. In all cases, the expected higher order convergence is observed. \textit{Right}: The error on the bulk domain is plotted on the y-axis, whilst now the x-axis shows the runtime of the computation. The higher order methods lead to shorter runtimes when high accuracy is demanded.}
   \label{fig:conv_lin} \vspace{-0.4cm}
\end{figure}
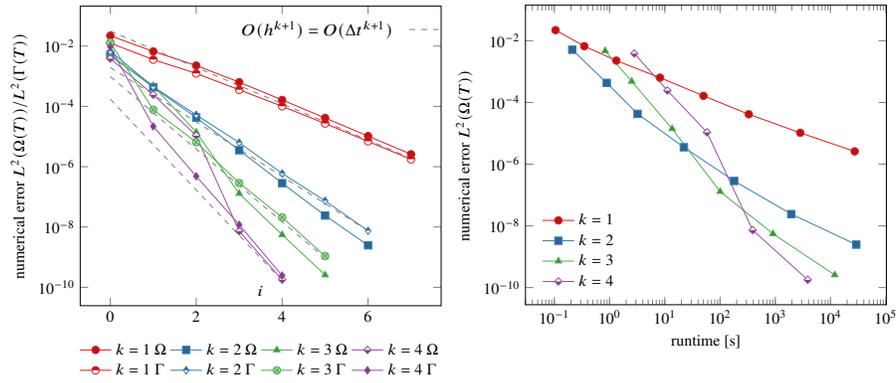

For the linear case ($b_{BS} = 0$), we show the convergence results in these norms in \Cref{fig:conv_lin}.
In all cases, we observe optimal order with discrete spaces of order $k$ in the $L^2$-norm, i.e.
\begin{align*}
\| u_B - u_{B,h} \|_{L^2(\Omega(T))}, ~ \| u_S - u_{S,h} \|_{L^2(\Omega(\Gamma))} \simeq \mathcal{O}(\Delta t^{k+1}) \simeq \mathcal{O}(h^{k+1}).
\end{align*}
In relation to the computation time relative to the resulting numerical error, as displayed in \Cref{fig:conv_lin}, right hand side, we note that the higher order methods show an error scaling behaviour benefitial for high accuracy demand simulations.\footnote{Note that we used the direct solver \texttt{umfpack} and a non-optimised parallisation for assembly. Hence, the specific timinig numbers should be rather interpreted as examples to show the scaling behaviour, and improvements might be possible.}

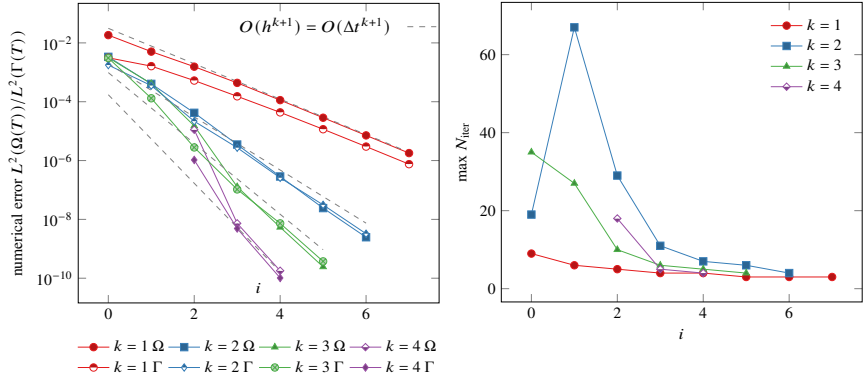
\begin{figure}
  \begin{tikzpicture}[scale=0.7]
    \begin{semilogyaxis}[ xlabel=$i$, ylabel={numerical error $L^2(\Omega(T))/ L^2(\Gamma(T))$}, legend entries ={ $k=1$ $\Omega$, $k=2$ $\Omega$, $k=3$ $\Omega$, $k=4$ $\Omega$, $k=1$ $\Gamma$, $k=2$ $\Gamma$, $k=3$ $\Gamma$, $k=4$ $\Gamma$
%      {[text width=2.5cm, text depth=] $O(h^{k+1})= O(\Delta t^{k+1})$}
    }, legend style={at={(0.5,-0.1)},anchor=north,legend columns=4, draw=none, legend cell align=left, fill=none},
    x label style={at={(axis description cs:0.5,0.05)},anchor=center},
%     ymax=3, ymin=3e-11
     ]
     \addplot table[x index = 0, y index =1] {out/conv_ks1_kt1_ktls1_from0_to7.dat};
     \addplot table[x index = 0, y index =1] {out/conv_ks2_kt2_ktls2_from0_to6.dat};
     \addplot table[x index = 0, y index =1] {out/conv_ks3_kt3_ktls3_from0_to5.dat};
     \addplot table[x index = 0, y index =1] {out/conv_ks4_kt4_ktls4_from2_to4.dat};
     \pgfplotsset{cycle list shift=1}
     \addplot table[x index = 0, y index =2] {out/conv_ks1_kt1_ktls1_from0_to7.dat};
     \addplot table[x index = 0, y index =2] {out/conv_ks2_kt2_ktls2_from0_to6.dat};
     \addplot table[x index = 0, y index =2] {out/conv_ks3_kt3_ktls3_from0_to5.dat};
     \addplot table[x index = 0, y index =2] {out/conv_ks4_kt4_ktls4_from2_to4.dat};
      \addplot[gray, dashed, domain=0:7] {(1/2^(x+2.5)))^2};
      \addplot[gray, dashed, domain=0:6] {(1/2^(x+3)))^3};
      \addplot[gray, dashed, domain=0:5] {(1/2^(x+2.5)))^4};
      \addplot[gray, dashed, domain=0:4] {(1/2^(x+2.5)))^5};
    \end{semilogyaxis}
     \node[scale=0.75] at (4.5,5.25) {$O(h^{k+1})= O(\Delta t^{k+1})$};
     \draw[scale=0.75, gray, dash=on 2.25pt off 2.25pt phase 0pt, line width=0.4*0.75pt] (6.25/0.75,5.25/0.75) -- (6.8/0.75,5.25/0.75);
   \end{tikzpicture}
     \begin{tikzpicture}[scale=0.7]
    \begin{axis}[ xlabel={$i$}, ylabel=max $N_{\text{iter}}$, legend entries ={ $k=1$, $k=2$, $k=3$, $k=4$
%      {[text width=2.5cm, text depth=] $O(h^{k+1})= O(\Delta t^{k+1})$}
    }, legend style={anchor=north,legend columns=1, draw=none, legend cell align=left, fill=none}, legend pos =north east,
%     x label style={at={(axis description cs:0.5,0.05)},anchor=center},
%     ymax=3, ymin=3e-11
     ]
     \addplot table[x index = 0, y index =3] {out/conv_ks1_kt1_ktls1_from0_to7.dat};
     \addplot table[x index = 0, y index =3] {out/conv_ks2_kt2_ktls2_from0_to6.dat};
     \addplot table[x index = 0, y index =3] {out/conv_ks3_kt3_ktls3_from0_to5.dat};
     \addplot table[x index = 0, y index =3] {out/conv_ks4_kt4_ktls4_from2_to4.dat};
%       \addplot[gray, dashed, domain=0:6] {(1/2^(x+2.5)))^2};
%       \addplot[gray, dashed, domain=0:5] {(1/2^(x+2.5)))^3};
%       \addplot[gray, dashed, domain=0:4] {(1/2^(x+2.5)))^4};
%       \addplot[gray, dashed, domain=0:3] {(1/2^(x+2.5)))^5};
    \end{axis}
    \draw[opacity=0] (0,0) -- (0,-1.5);
%      \node[scale=0.75] at (4.5,5.25) {$O(h^{k+1})= O(\Delta t^{k+1})$};
%      \draw[scale=0.75, gray, dash=on 2.25pt off 2.25pt phase 0pt, line width=0.4*0.75pt] (6.25/0.75,5.25/0.75) -- (6.8/0.75,5.25/0.75);
   \end{tikzpicture}
   \caption{Numerical convergence results for the non-linear model system ($b_{BS} = 1$), where the nonlinearity is solved by a Newton solver. \textit{Left}: Numerical errors in $L^2$ norms at final time on bulk and surface. \textit{Right}: Numbers of Newton iterations required. For coarse meshes/ large timesteps, higher iteration numbers might be needed for higher order.}
   \label{fig:conv_nonlin}\vspace{-0.4cm}
\end{figure}

We also observe the same higher order convergence for the non-linear problem, $b_{BS} = 1$, c.f. \Cref{fig:conv_nonlin} left-hand side.
Also, in the right-hand side of this figure, we evaluate the required number of Newton iterations and consider the maximal number over all time steps. We observe that for coarse meshes/ large time steps and higher order methods, the number of required Newton iterations might become high. For $k=4$ and $i=0,1$, no convergence in Newton was observed. In the future, one could try to improve this pre-asymptotic behaviour of the proposed method, e.g. by damped Newton methods or alike.
\section{Summary and Outlook} \label{sec:summ_outlook}
To summarise, we have presented an unfitted space-time finite element method for a coupled surface-bulk problem. Overall, the techniques developed for the bulk and surface problems, respectively, transfer. In particular, for the representation of the discrete geometry with higher order accuracy, the isoparametric mapping of \cite{HLP2022, L_CMAME_2016, heimann2023geometrically} can be applied here as well, yielding a higher order accurate geometry for arbitrary order parameters in space and time, together with suiting numerical integration routines. For the time discretisation, we opted for tensor-product space-time finite elements, which allow for a flexible choice of discretisation order in space and time and a time stepping procedure. Beyond time slice boundaries, continuity is enforced weakly with an upwind-style penalty term. For stabilisation, we combine the space-time direct variant of the Ghost penalty stabilisation \cite{burman2010ghost, burman15, preuss18, heimann20} on the bulk with a surface normal gradient stabilisation on the ambient elements\cite{grande2018analysis, burman18cutfemcodim, sassreusken23, reusken2024analysis}. Overall, the resulting method shows convergence of order $k+1$ for polynomial order choice $k$ for the discrete spaces and geometrical accuracy, when the error is measured in an $L^2$ norm at final time on bulk or surface.

In the future, one could aim at performing a numerical analysis showing these convergence properties in a rigorous mathematical manner. (See \cite{preuss18, heimann20} for related bulk, and \cite{reusken2024analysis} for surface results) In relation to the non-linear problem, the suggested Newton's method might allow for further improvements to yield faster convergence in the context of coarse meshes/ time steps. As far as computations are concerned, the suggested methods might also be transferred to other physical problems than the considered convection-diffusion model problem on bulk and surface.
\subsection*{Acknowledgements}
The author acknowledges helpful suggestions by Christoph Lehrenfeld on this text and discussions on the topic.
\bibliographystyle{abbrv}
\bibliography{paper}
\end{document}